\documentclass{article}
\usepackage{amssymb}
\usepackage[T1]{fontenc}

\begin{document}
\newcommand{\ov}{\overline}
\newcommand{\un}{\underline}
\newcommand{\proof}{\bf {Proof:} \rm}
\newcommand{\dirac}{ D_{h}}
\newcommand{\dif}{ \partial_{h}}
\newcommand{\lapl}{ \Delta_{h}}
\newcommand{\euler}{ E_{h}^{\pm}}
\newcommand{\gam}{ \Gamma_{h}^{\pm}}

\newcommand{\BR}{{\mathbb R}}
\newcommand{\BN}{{\mathbb N}}
\newcommand{\BZ}{{\mathbb Z}}
\newcommand{\BH}{{\mathbb H}}

\newcommand{\Grad}{{\rm grad}}
\newcommand{\Curl}{{\rm curl}}
\newcommand{\Div}{{\rm div}}
\newcommand{\Sc}{{\rm Sc\,}}
\newcommand{\Vec}{{\rm Vec\,}}
\newcommand{\wzwo} {\stackrel{\circ}{{\cal W}^1_2}}
\newcommand{\wzwoz} {\stackrel{\circ}{{\cal W}^2_2}}
\newcommand{\wzwok}{\stackrel{\circ}{\cW_2^k}}
\newcommand{\wzwokl}{\stackrel{\circ}{\cW_2^{k,l}}}
\newcommand{\G}{\mathcal{G}}

\newcommand{\e}{{\bf e}}
\newcommand{\cW}{{W}}
\newcommand{\cL}{{L}}
\newcommand{\cC}{{C}}
\newcommand{\cM}{{M}}

\newcommand{\cl}{C \kern -0.1em \ell}

\newcommand{\qed}{$\blacksquare$}
\newtheorem{theorem}{Theorem}[section]
\newtheorem{remark}{Remark}[section]
\newtheorem{lemma}{Lemma}[section]
\newtheorem{proposition}{Proposition}[section]
\newtheorem{corollary}{Corollary}[section]
\newtheorem{definition}{Definition}[section]
\newtheorem{example}{Example}[section]

\title{Fischer Decomposition for Difference Dirac Operators}
\author{N. Faustino\thanks {supported by {\it Funda\c c\~ao para Ci\^encia e Tecnologia} under PhD-grant no. SFRH/BD/17657/2004} \\ Departamento de Matem\'{a}tica,\\
Universidade de Aveiro,\\
P-3810-193 Aveiro,\\ Portugal \and U. K\"ahler\thanks {Research
(partially) supported by {\it Unidade de
Investiga\c c\~ao Matem\'atica e Aplica\c c\~oes} of the University of Aveiro.}\\
Departamento de Matem\'{a}tica,\\
Universidade de Aveiro,\\
P-3810-193 Aveiro,\\ Portugal} \maketitle

\begin{quote}
\begin{abstract}
We establish the basis of a discrete function theory starting with a
Fischer decomposition for difference Dirac operators. Discrete
versions of homogeneous polynomials, Euler and Gamma operators are
obtained. As a consequence we obtain a Fischer decomposition for the
discrete Laplacian.
\end{abstract}
\end{quote}


{\bf keywords}: Difference Dirac operator, Difference Euler
operator, Fischer decomposition

\section{Introduction}

Clifford analysis is a powerful tool to solve some kinds of problems related with vector field analysis. \\
A comprehensive description of Clifford function theory was
introduced by F. Brackx, R. Delanghe and F. Sommen in \cite{Brackx}
and later by R. Delanghe, F. Sommen and V. Sou\u{c}ek
in~\cite{Delanghe}.

In~\cite{GS,GS2}, K. G\"urlebeck and W. Spr\"o\ss ig proposed some
strategies to solve boundary value problems based on the study of
existence, uniqueness, representation, and regularity of solutions
with the help of an operator calculus. In the same books, the
authors introduce also the basic ideas to develop a discrete
counterpart of the continuous treatment of boundary value problems
with the introduction of a discrete operator calculus in order to
find a well-adapted numerical approach. An explicit discrete version
of the Borel-Pompeiu formula was presented for dimension $n=3$.

This was further developed in~\cite{GH,GH2}, where K. G\"urlebeck
and A. Hommel developed finite difference potential methods in
lattice domains based on the concept of discrete fundamental
solutions for the difference Dirac operator which generalizes the work
developed by Ryabenkij in \cite{Ryaij}. A numerical application of
this theory was presented recently by N. Faustino, K. G\"urlebeck,
A. Hommel, and U. K\"ahler in~\cite{Faust} for the incompressible
stationary Navier-Stokes equations. In this paper, the authors
proposed a scheme which solves efficiently problems in unbounded
domains and show the convergence of the numerical scheme for
functions with H\"older regularity which is a better gain compared
with the convergence results for classical difference schemes.

Moreover, while all these papers claim to be based on discrete
function theoretical approaches, from the concepts of the theory of
monogenic functions only the Borel-Pompeiu formula and with it
Cauchy's integral formula were obtained. There is no ``real''
development of a discrete monogenic function theory up to now.

This paper is supposed to be a step in this direction. To this end
discrete versions of a Fischer decomposition, Euler and Gamma
operators are obtained. For the sake of simplicity we consider in
the first part only Dirac operators which contain only forward or
backward finite differences. Of course, these Dirac operators do not
factorize the classic discrete Laplacian. Therefore, we will
consider in the last chapter a different definition of a difference
Dirac operator in the quaternionic case (c.f.~\cite{GH}) which do
factorizes the discrete Laplacian.

Let us emphasize in the end a major obstacle in the discrete case.
While in the continuous case the are only one partial derivative for
each coordinate $x_j$ we have two finite differences in the discrete
case. Therefore, we will have not only one Euler or Gamma operator
as in the continuous case, but several. Each one will turn out to be
connected to one particular Dirac operator.

\section{Preliminaries}

Let $\e_1,\ldots,\e_n$ be an orthonormal basis of $\BR^n$. The
Clifford algebra $\cl_{0,n}$ is the free algebra over $\BR^n$
generated modulo the relation
$$
x^2=-|x|^2\e_0,
$$
where $\e_0$ is the identity of $\cl_{0,n}$. For the algebra
$\cl_{0,n}$ we have the anti-commutation relationship
$$
\e_i\e_j+\e_j\e_i=-2\delta_{ij}\e_0,
$$
where $\delta_{ij}$ is the Kronecker symbol. In the following we
will identify the Euclidean space $\BR^n$ with
$\bigwedge^1\cl_{0,n}$, the space of all vectors of $\cl_{0,n}$.
This means that each element $x$ of $\BR^n$ may be represented by
$$
x=\sum_{i=1}^n x_i\e_i.
$$
From an analysis viewpoint one extremely crucial property of the
algebra $\cl_{0,n}$ is that each non-zero vector $x\in\BR^n$ has a
multiplicative inverse given by $\frac{-x}{|x|^2}$. Up to a sign
this inverse corresponds to the Kelvin inverse of a vector in
Euclidean space. Moreover, given a general Clifford number
$a=\sum_A\e_Aa_A, A\subset\{1,\ldots,n\}$ we denote by $\Sc
a=a_\emptyset$ the scalar part and by $\Vec
a=\e_1a_1+\ldots+\e_na_n$ the vector part.

For all what follows let $\Omega\subset\BR^n$ be a bounded domain
with a sufficiently smooth boundary $\Gamma=\partial\Omega$. Then
any function $f:\Omega\mapsto \cl_{0,n}$ has a representation
$f=\sum_A\e_Af_A$ with $\BR$-valued components $f_A$. We now
introduce the Dirac operator $D=\sum_{i=1}^n\e_i\frac{\partial}
{\partial x_i}$. This operator is a hypercomplex analogue to the
complex Cauchy-Riemann operator. In particular we have that
$D^2=-\Delta$, where $\Delta$ is the Laplacian over $\BR^n$. A
function $f:\Omega\mapsto\cl_{0,n}$ is said to be {\it
left-monogenic} if it satisfies the equation $(Df)(x)=0$ for each
$x\in\Omega$. A similar definition can be given for right-monogenic
functions. Basic properties of the Dirac operator and left-monogenic
functions can be found
in~\cite{Brackx},~\cite{Delanghe},~\cite{GS2}, and~\cite{GS}.

Now, we need some more facts for our discrete setting. To discretize
point-wise the partial derivatives $\frac{\partial} {\partial x_i}$
in the equidistant lattice with mesh width $h>0$,
$\BR^n_h=\{mh=(m_1h,\ldots,m_nh):m\in\BZ^n\}$, we introduce
forward/backward differences $\dif^{\pm i}$:
\begin{equation}
\dif^{\pm i} u (mh)=\mp \frac{u(mh)-u(mh \pm h \e_i)}{h}
\end{equation}

These forward/backward differences $\dif^{\pm i}$ satisfy the
following product rules
\begin{eqnarray} \label{difRule1} (\dif^{\pm
i}fg)(mh)&=&f(mh)(\dif^{\pm i}g)(mh)+(\dif^{\pm i}f)(mh)g(mh \pm h\e_i),  \\
\label{difRule2}(\dif^{\pm i}fg)(mh)&=&f(mh \pm h\e_i)(\dif^{\pm
i}g)(mh)+(\dif^{\pm i}f)(mh)g(mh).
\end{eqnarray}

The forward/backward discretizations of the Dirac operator are given
by
\begin{equation}
\label{difD} \dirac^{\pm} =\sum_{i=1}^n \e_i \dif^{\pm i}.
\end{equation}

In this paper we will also use the following multi-index
abbreviations:
\begin{eqnarray*}
(mh)^{(\alpha)}:=(m_1h)^{\alpha_1}(m_2h)^{\alpha_2} \ldots
(m_nh)^{\alpha_n}; \\
\alpha ! := \alpha_1 ! \alpha_2 ! \ldots \alpha_n !; \\
|\alpha|:= \alpha_1+\alpha_2+ \ldots \alpha_n \\
\dif^{\pm \e_i}:=\dif^{\pm i}; \\
\dif^{\pm \alpha_i \e_i}:=\left(\dif^{\pm \e_i}\right)^{\alpha_i},
\end{eqnarray*}
for
$\alpha=(\alpha_1,\alpha_2,\ldots,\alpha_n)=\sum_{i=1}^n\e_i\alpha_i$.

\section{Fischer decomposition}
The basic idea of a Fischer decomposition is to decompose any
homogeneous polynomial into monogenic homogeneous polynomials of
lower degrees. In the classic case such a decomposition is based
on the fact that the powers $x^s$ are homogeneous and that
$\frac{\partial x_i^s}{\partial x_i}=sx_i^{s-1}$. A first idea
would be to consider instead of $x^s$ simply the powers $(mh)^s$,
but while these powers are still homogeneous the last condition is
not true in the discrete case, unfortunately. Therefore, we will
start by introducing discrete homogeneous powers which will play
the equivalent role of $x^s$ in the discrete case.

\subsection{Multi-index factorial powers} \label{homogeneous}

Starting from the one-dimensional factorial powers
\begin{eqnarray}\label{fact}(m_ih)^{(0)}_{\mp}:=1,(m_ih)^{(s)}_{\mp}:=\prod_{k=0}^{s-1}(m_ih
\mp kh), s \in \BN\end{eqnarray} we introduce the multi-index
factorial powers of degree $|\alpha|$ by
$$ (mh)^{(\alpha)}_{\mp}=\prod_{i=1}^n (m_ih)^{(\alpha_i)}_{\mp}.
$$

The one-dimensional factorial powers $(m_ih)^{(s)}_{\mp}$ have the
following properties
\begin{itemize}
\item[\bf P1.] $(m_ih)^{(s+1)}_{\mp}=(m_ih \mp sh)(m_ih)^{(s)}_{\mp}$;
\item[\bf P2.] $\dif^{\pm
j}(m_ih)^{(s)}_{\mp}=s(m_ih)^{(s-1)}_{\mp}\delta_{i,j}$;
\item[\bf P3.] $\dif^{\mp
j}(m_ih)^{(s)}_{\mp}=s(m_ih \mp h)^{(s-1)}_{\mp}\delta_{i,j}$;
\item[\bf P4.] $(m_ih)^{(s)}_{\mp} \rightarrow x_i^s=(m_ih)^s$ for $h \rightarrow
0$,
\end{itemize}
where $\delta_{i,j}$ denotes the standard Kronecker symbol.

As a direct consequence of these properties, we obtain the following lemmas:
\begin{lemma} \label{prop1}
The multi-index factorial powers of degree $|\alpha|$,
$(mh)^{(\alpha)}_{\mp}$, satisfy
$$\sum_{i=1}^n (m_ih) \dif^{\pm i} (mh \mp h \e_i)^{(\alpha)}_{\mp}=|\alpha|(mh)^{(\alpha)}_{\mp}.$$
\end{lemma}

\begin{lemma} \label{prop2}
The multi-index factorial powers of degree $|\alpha|$,
$(mh)^{(\alpha)}_{\mp}$, satisfy
$$ \dif^{\pm \beta} (mh)^{(\alpha)}_{\mp}=\alpha ! \delta_{\alpha,\beta} \qquad \mbox{for }|\beta|=|\alpha|.$$
\end{lemma}

\begin{lemma} \label{prop3} The multi-index factorial powers $(mh)^{(\alpha)}_{\mp}$ of
degree $|\alpha|$ approximate the classical multi-index powers
$x^{(\alpha)}$ of degree $|\alpha|$ at each point $x=mh$.\end{lemma}

Let us remark that we have the following relationships between the
multi-index factorial powers and the usual powers:

\begin{theorem} \label{comb1}
The powers $(m_ih)^{\alpha_i}$ and $(m_ih)_{\mp}^{\alpha_i}$ are
related by
$$(m_ih)^{\alpha_i}_{\mp}=
\sum_{k_i=0}^{\alpha_i}S^{\alpha_i}_{k_i}(m_ih)^{k_i}$$
$$(m_ih)^{\alpha_i}= \sum_{k_i=0}^{\alpha_i}T^{\alpha_i}_{k_i}(m_ih)_{\mp}^{k_i},$$
where $S^{\alpha_i}_{k_i}$ are the Stirling numbers of the first
kind and $T^{\alpha_i}_{k_i}$ are the Stirling numbers of the second
kind.
\end{theorem}

The sketch of the proof of this theorem can be found, e.g., in
\cite{Laks}.

\begin{theorem} \label{comb2}
The multi-index powers $(mh)^{(\alpha)}$ and $(mh)^{(\alpha)}_{\mp}$
are related by
\begin{eqnarray}\label{Stirl1}(mh)^{(\alpha)}_{\mp}=
\sum_{|\beta|=0}^{|\alpha|}K_{\beta}^{\alpha}(mh)^{(\beta)}, \\
\label{Stirl2}(mh)^{(\alpha)}=
\sum_{|\beta|=0}^{|\alpha|}L_{\beta}^{\alpha}(mh)^{(\beta)}_{\mp}.
\end{eqnarray}
Moreover,
$$K_{\beta}^{\alpha}=\sum_{l_{n-1}=0}^{|\beta|}
\sum_{l_{n-2}=0}^{l_{n-1}} \ldots
\sum_{l_{1}=0}^{l_2}S^{\alpha_1}_{l_1}S^{\alpha_2}_{l_2-l_1}
\ldots S^{\alpha_n}_{|\beta|-l_{n-1}} ,$$
$$L_{\beta}^{\alpha}=\sum_{l_{n-1}=0}^{|\beta|}
\sum_{l_{n-2}=0}^{l_{n-1}} \ldots
\sum_{l_{1}=0}^{l_2}T^{\alpha_1}_{l_1}T^{\alpha_2}_{l_2-l_1}
\ldots T^{\alpha_n}_{|\beta|-l_{n-1}} .$$
\end{theorem}

We will just prove identity (\ref{Stirl1}). The proof of identity
(\ref{Stirl2}) is analogous to the proof of identity (\ref{Stirl1}).

\proof Using Theorem~\ref{comb1} and multiplying the polynomials
$(m_1h)^{\alpha_1}$ and $(m_2h)^{\alpha_2}$, we obtain
$$ (m_1h)^{\alpha_1}(m_2h)^{\alpha_2}=\sum_{\beta_1+\beta_2=0}^{\alpha_1+\alpha_2}
K_{\alpha_1,\alpha_2}^{\beta_1,\beta_2}(m_1h)_{\mp}^{(\beta_1)}(m_2h)_{\mp}^{(\beta_2)}$$
with
$K_{\alpha_1,\alpha_2}^{\beta_1,\beta_2}=\left(\sum_{l_1=0}^{\beta_1+\beta_2}
S^{\alpha_1}_{l_1}S^{\alpha_2}_{\beta_1+\beta_2-l_1} \right)$.

Using again Theorem~\ref{comb1} and multiplying the polynomials
$(m_1h)^{\alpha_1}(m_2h)^{\alpha_2}$ and $(m_3h)^{\alpha_3}$, we
obtain
$$ (m_1h)^{\alpha_1}(m_2h)^{\alpha_2}(m_3h)^{\alpha_3}=\sum_{\beta_1+\beta_2+\beta_3=0}^{\alpha_1+\alpha_2+\alpha_3}
K_{\alpha_1,\alpha_2,\alpha_3}^{\beta_1,\beta_2,\beta_3}(m_1h)_{\mp}^{(\beta_1)}(m_2h)_{\mp}^{(\beta_2)}(m_3h)_{\mp}^{(\beta_3)}$$
with
$K_{\alpha_1,\alpha_2,\alpha_3}^{\beta_1,\beta_2,\beta_3}=\sum_{l_2=0}^{\beta_1+\beta_2+\beta_3}
\sum_{l_1=0}^{l_2}S^{\alpha_1}_{l_1}S^{\alpha_2}_{l_2-l_1}S^{\alpha_3}_{\beta_1+\beta_2+\beta_3-l_2}$.

Applying this procedure recursively, we obtain
\begin{eqnarray}
\label{multiexp}(mh)^{(\alpha)}&=&\sum_{|\beta|=0}^{|\alpha|}K_{\beta}^{\alpha}(mh)_{\mp}^{(\alpha)}
\end{eqnarray}
with $K_{\beta}^{\alpha}=\sum_{l_{n-1}=0}^{|\beta|}
\sum_{l_{n-2}=0}^{l_{n-1}} \ldots
\sum_{l_{1}=0}^{l_2}S^{\alpha_1}_{l_1}S^{\alpha_2}_{l_2-l_1}
\ldots S^{\alpha_n}_{|\beta|-l_{n-1}}.$ \qed

For all what follows, let $\Pi_{k}^{\pm}$ denote the space of all
Clifford-valued polynomials of degree $k$, $P^{\pm}_k$, generated by
the powers $(mh)^{(\alpha)}_{\mp}$ of degree $|\alpha|=k$, and
$\Pi^{\pm}$ be the countable union of all Clifford-valued
polynomials of degree $k \geq 0$. Furthermore, let ${\cal
M}_k^{\pm}=\Pi^{\pm}_k \cap \ker \dirac^{\pm}$ be the space of
discrete monogenic polynomials of degree $k$. Based on
Lemma~\ref{prop1},~\ref{prop2} and~\ref{prop3}, we will show that it
is possible to obtain discrete versions for the Fischer
decomposition as well as define discrete versions of the Euler and
Gamma operators.

\subsection{The main theorem}

For two Clifford-valued polynomials of degree $k$, $P^{\pm}_k$ and
$Q^{\pm}_k \in\Pi^{\pm}_k$ given by
\begin{eqnarray*}
P^{\pm}_k(mh) & = & \sum_{|\alpha|=k}(mh)_{\mp}^{(\alpha)}a_{\alpha}^{\pm}\\
Q^{\pm}_k(mh) & = & \sum_{|\alpha|=k}(mh)_{\mp}^{(\alpha)}b_{\alpha}^{\pm}\\
\end{eqnarray*}
we define the Fischer inner product by
\begin{equation}\label{fischer}[P^{\pm}_k,Q^{\pm}_k]_{h}:=\sum_{|\alpha|=k} \alpha
! \Sc (\overline{a_{\alpha}^{\pm}}b_{\alpha}^{\pm}).
\end{equation}

Denoting by $P^{\pm}_k(\dirac^{\pm})$ the difference operator
obtained from the polynomial $P^{\pm}_k$ in powers of $mh$ by
replacing $m_ih$ by the difference operator $\dif^{\pm i}$ (just
like in the continuous case, c.f.~\cite{Delanghe}), we have by
Lemma~\ref{prop2} the identity
\begin{equation}
\label{innerfischer}[P^{\pm}_k,Q^{\pm}_k]_{h}:=\Sc(\overline{P^{\pm}_k(\dirac^{\pm})}Q^{\pm}_k)(0)
\qquad P^{\pm}_k,Q^{\pm}_k\in\Pi^{\pm}_k.
\end{equation}
With other words, we can express the Fischer inner product by
applying the difference operator $P^{\pm}_k(\dirac^{\pm})$ to the
polynomial $P^{\pm}_k$ and evaluate the scalar part at the point
$mh=0$.

Moreover, due to $\overline{\dirac^{\pm}}=- \dirac^{\pm}$ the
Fischer inner product has the important property:
\begin{equation}\label{fischer3}
[(mh) P^{\pm}_k,Q^{\pm}_k]_{h}=-[P^{\pm}_k,\dirac^{\pm}
Q^{\pm}_k]_h.
\end{equation}

This property allows us to prove the following theorem:
\begin{theorem}\label{fisher2}
We have
$$
\Pi_k^{\pm}={\cal M}_k^{\pm}+(mh) \Pi_{k-1}^{\pm}.
$$
Moreover, the subspaces ${\cal M}_k^{\pm}$ and $(mh) \Pi_{k-1}^{\pm}
$ are orthogonal with respect to the Fischer inner product.
\end{theorem}

Before we prove Theorem~\ref{fisher2}, we will prove the following
inclusion property:

\begin{lemma} \label{prop5}
For the set $\dirac^{\pm} \Pi_k^{\pm}:=\{\dirac^{\pm} P_k^{\pm}:
P_k^{\pm} \in \Pi_k^{\pm} \}$, we have the inclusion:
$$ \dirac^{\pm} \Pi_k^{\pm}:=\{\dirac^{\pm} P_k^{\pm}: P_k^{\pm} \in
\Pi_k^{\pm} \} \subset \Pi_{k-1}^{\pm}.
$$
\end{lemma}

\proof

Let
$P^{\pm}_k(mh)=\sum_{|\alpha|=k}(mh)_{\mp}^{(\alpha)}a_{\alpha}^{\pm}
\in \Pi_k^{\pm}$. Applying $\dirac^{\pm}$, we obtain from $\dif^{\pm
i} (mh)_{\mp}^{(\alpha)}=\alpha_i(mh)_{\mp}^{(\alpha-\e_i)}$, the
identity
\begin{equation}\label{DP}(\dirac^{\pm} P^{\pm}_k)(mh)=\sum_{i=1}^n
\sum_{|\alpha|=k}(mh)_{\mp}^{(\alpha-\e_i)} \alpha_i\e_i
a_{\alpha}^{\pm}
\end{equation}
Because $\alpha_i\e_i a_{\alpha}^{\pm}$ is a Clifford constant we
have a linear combination of polynomials of degree
$|\alpha-\e_i|=k-1$ on the right hand side of (\ref{DP}). Hence,
$\dirac^{\pm} P^{\pm}_k \in \Pi^{\pm}_{k-1}$. \qed

\proof({\bf Theorem \ref{fisher2}}) Because of $\Pi_k^{\pm}= (mh)
\Pi_{k-1}^{\pm}+\left((mh ) \Pi_{k-1}^{\pm} \right)^\bot$ it is
enough to prove that $\left((mh ) \Pi_{k-1}^{\pm}\right)^\bot={\cal
M}_{k-1}^{\pm}$. For this we choose $P_{k-1}^{\pm}\in
\Pi_{k-1}^{\pm}$ arbitrarily and assume that for some
$P_k^{\pm}\in\Pi_k^{\pm}$ we have
$$
[(mh) P_{k-1}^{\pm},P_k^{\pm}]_h=0.
$$
Due to (\ref{fischer}) we have $[P_{k-1}^{\pm},\dirac^{\pm}
P_k^{\pm}]_h=0$ for all $P_{k-1}^{\pm}$. As $\dirac^{\pm}
P_k^{\pm}\in\Pi_{k-1}^{\pm}$ by Lemma~\ref{prop5}, we obtain
$\dirac^{\pm} P_k^{\pm}=0$ or $P_k^{\pm}\in {\cal M}_k^{\pm}$. This
means that $((mh) \Pi_{k-1}^{\pm})^\bot\subset {\cal M}_{k}^{\pm}$.
Now, let $P_k^{\pm}\in {\cal M}_k^{\pm}$. Then we have for each
$P_{k-1}^{\pm}\in\Pi_{k-1}^{\pm}$
\begin{eqnarray*}
[(mh)
P_{k-1}^{\pm},P_k^{\pm}]_h & = & -[P_{k-1}^{\pm},\dirac^{\pm}P_k^{\pm}]_h\\
& = & 0
\end{eqnarray*}
and, therefore, $((mh) \Pi_{k-1}^{\pm})^\top={\cal M}_{k-1}^{\pm}$.
\qed

From this theorem we obtain the Fischer decomposition with respect
to our difference Dirac operators $\dirac^{\pm}$.
\begin{theorem}{\bf Fischer decomposition} Let $P^{\pm}_k \in \Pi_{k}^{\pm}$  then
\begin{eqnarray}
P^{\pm}_k(mh)&=&\sum_{s=0}^{k-1}(mh)^sM_{k-s}^{\pm}(mh).
\end{eqnarray}
where each $M^{\pm}_j$ denotes the homogeneous discrete monogenic
polynomials of degree $j$ with respect to the Dirac operators
$\dirac^{\pm}$.
\end{theorem}

\subsection{Difference Euler and Gamma operators} \label{EulerGamma}

Based on Lemma~\ref{prop1} we will introduce discrete versions of
the Euler and Gamma operators presented in~\cite{Delanghe}.

First of all, we introduce the second order difference operator
$A_h^{\pm}$ by
\begin{eqnarray}
\label{Ah}A_h^{\pm} = \mp h \sum_{i=1}^n (m_ih) \dif^{\pm
i}\dif^{\mp i}.
\end{eqnarray}

\begin{definition}
For a lattice function $f_h:\BR^n_h \rightarrow \cl_{0,n}$, we
introduce the difference Euler operator $\euler$ by
$$ (\euler f_h)(mh)= \sum_{i=1}^n (m_ih)  (\dif^{\pm i} f_h)(mh \mp h \e_i)$$
and the difference Gamma operator $\gam$ by
$$(\gam f_h)(mh)=-\sum_{j<k}\e_j\e_k (L_{jk}^{\pm}f_h)(mh) - (A_h^{\pm} f_h)(mh),$$
where $L_{jk}^{\pm}:=(m_jh)\dif^{\pm k}-(m_kh) \dif^{\pm j}$.
\end{definition}

It looks surprising that we have in the definition of the Gamma
operator a term which contains second order differences, but we
would like to remark that for $h\rightarrow 0$ this term vanishes
and we will get the usual continuous Gamma operator. As a matter of
fact this term arises due to the fact that in the discrete case
translations are involved in the definition of finite differences/
finite difference operators.

Using the definition of the difference Euler operator and
Lemma~\ref{prop1}, we obtain for polynomials homogeneous of degree
$k$, $P^{\pm}_{k} \in \Pi^{\pm}_{k}$, $\euler
P^{\pm}_{k}=kP^{\pm}_{k}$, and, moreover, we can show that a
function $f_h$ homogeneous of degree $k$ satisfy $\euler f_h=k f_h$.
This fact provides a good motivation for calling $\euler$ Euler
operator, i.e. an operator who measures the degree of homogeneity of
a homogeneous function.

It follows from the definition of the Euler and Gamma operator that
\begin{eqnarray}
\label{id1}(mh)\dirac^{\pm}f_h&=&-\sum_{i=1}^n (m_ih) \dif^{\pm
i}f_h+\sum_{j<k}\e_j\e_k L_{jk}^{\pm}f_h \\
\label{id2}&=&-(\euler+\gam )f_h
\end{eqnarray}
The proof of (\ref{id2}) is easily obtained by adding and
subtracting the term
$$ -  \sum_{i=1}^n  (m_ih) \dif^{\pm i}(f_h(mh)-f_h(mh \mp h \e_i)),$$
which is the same as $A_h^{\pm} f_h$, on the right hand side of
(\ref{id1}). Moreover, for discrete monogenic polynomials of degree
$k$, $M_k^{\pm} \in {\cal M}_k^{\pm}$, we have $\gam
M_k^{\pm}=-kM_k^{\pm}$.

For all what follows, we introduce the difference operators
\begin{eqnarray}
\label{Bh}B_h^{\pm}=\pm h \sum_{i=1}^n \dif^{\pm i}, \\
\label{Ch}C_h^{\pm}f_h=\sum_{i=1}^n (m_ih)\e_if_h(\cdot \mp h \e_i) \\
\label{Rh}R_{h,r}^{\pm}=rI+\euler-A_h^{\pm}, \\
\label{Vh}V_{h,r}^{\pm}=R_{h,r}^{\pm}+\frac{1}{2}B_h^{\pm}.
\end{eqnarray}
where $I$ is the identity operator and $r$ a real number.

From the identity
\begin{eqnarray*} \left(
(mh)\dirac^{\pm}+\dirac^{\pm} (mh) \right)f_h&=&-2\sum_{j=1}^n
(m_jh) \dif^{\pm
j}f_h-nf_h\nonumber \\
&=&-2(\euler-A_h^{\pm})f_h-nf_h \\
&=&-2R_{h,n/2}^{\pm}f_h
\end{eqnarray*}
we get \begin{eqnarray}\label{Dmh} \left(\dirac^{\pm}
(mh)\right)f_h=(-2R_{h,n/2}^{\pm}+\euler+\gam) f_h,\end{eqnarray}
by applying identity (\ref{id2}).

Now, we will show some important facts regarding our difference operators.
\begin{proposition} \label{diracE}
For a lattice function $f_h:\BR^n_h \rightarrow \cl_{0,n}$, we have
$$ \dirac^{\pm}\euler f_h=\dirac^{\pm} f_h+\euler \dirac^{\pm} f_h .$$
\end{proposition}

\proof

Starting from the definition, we can split $\dirac^{\pm}\euler
f_h$ in the sum
\begin{eqnarray}
(\dirac^{\pm}\euler f_h)(mh)&=& I_1^{\pm}(mh)+I_2^{\pm}(mh)
\end{eqnarray}
with
$$I_1^{\pm}(mh)=\sum_{i=1}^n \e_i\dif^{\pm i}\left((m_ih)  (\dif^{\pm i}f_h)(mh \mp h \e_i) \right)$$
and
$$I_2^{\pm}(mh)=\sum_{j=1}^n\sum_{i \neq j} \e_j\dif^{\pm j}\left( (m_ih) (\dif^{\pm i}f_h)(mh \mp h \e_i)\right).$$
Applying the product rule for finite differences (\ref{difRule1}) in $I^{\pm}_1$, we
obtain
\begin{eqnarray*}I_1^{\pm}(mh)&=& \sum_{i=1}^n \e_i\left((\dif^{\pm i}f_h)(mh)+(m_ih) (\dif^{\pm 2
\e_i}f_h)(mh \mp h\e_i) \right)\\
&=&(\dirac^{\pm}f_h)(mh)+\sum_{i=1}^n \e_i(m_ih) (\dif^{\pm 2
\e_i}f_h)(mh \mp h\e_i).
\end{eqnarray*}
On the other hand, $$I^{\pm}_2(mh)=\sum_{j=1}^n\e_j\sum_{i \neq
j}(m_ih) \left( (\dif^{\pm i}\dif^{\pm j}f_h)(mh \mp h
\e_i)\right).$$ Thus, we have
\begin{eqnarray*}
(\dirac^{\pm}\euler f_h)(mh)&=&(\dirac^{\pm}f_h)(mh)+
\sum_{i,j=1}^n\e_j(m_ih)  (\dif^{\pm i}\dif^{\pm j}f_h)(mh \mp h
\e_i)
\\
&=&(\dirac^{\pm}f_h)(mh)+(\euler\dirac^{\pm} f_h)(mh).
\end{eqnarray*}
\qed

\begin{proposition} \label{diracmh}
For a lattice function $f_h:\BR^n_h \rightarrow \cl_{0,n}$, we have
\begin{equation}\label{DE}
\dirac^{\pm}((mh)f_h)=-2V^{\pm}_{h,n/2}f_h-(mh)\dirac^{\pm}f_h\end{equation}
\end{proposition}

\proof Using the product rule for finite differences
(\ref{difRule1}) and the identity $-2m_ih=\e_i(mh)+(mh)\e_i,
i=1,\ldots,n,$ we get
\begin{eqnarray*}
\dirac^{\pm}((mh)f_h(mh))&=&-\sum_{i=1}^n f_h(mh \pm h
\e_i)-2\sum_{i=1}^n (m_ih)\dif^{\pm i}f_h(mh)\\
&&-(mh)(\dirac^{\pm}f_h)(mh) \\
&=&-nf_h(mh)-\sum_{i=1}^n (f_h(mh \pm h \e_i)-f_h(mh))- \\
&&-2(\euler f_h- A_h^{\pm}f_h)(mh)-(mh)(\dirac^{\pm}f_h)(mh) \\
&=&-nf_h(mh)-2(\euler f_h- A_h^{\pm}f_h)(mh)- \\ &&-(B_h^{\pm}f_h)(mh)-(mh)(\dirac^{\pm}f_h)(mh) \\
&=&-2(R^{\pm}_{h,n/2}f_h+\frac{1}{2}B_h^{\pm}f_h)(mh)-(mh)(\dirac^{\pm}f_h)(mh) \\
&=&-2(V^{\pm}_{h,n/2}f_h)(mh)-(mh)(\dirac^{\pm}f_h)(mh).
\end{eqnarray*}
 \qed

From Proposition~\ref{diracmh} and from the commutation properties
$\dirac^{\pm}A_h^{\pm}=A_h^{\pm}\dirac^{\pm}$ and
$\dirac^{\pm}B_h^{\pm}=B_h^{\pm}\dirac^{\pm}$ follow the operator
relations
\begin{eqnarray}
\label{DR}\dirac^{\pm}R_{h,r}^{\pm}=R_{h,r+1}^{\pm}\dirac^{\pm}, \\
\label{DV}\dirac^{\pm}V_{h,r}^{\pm}=V_{h,r+1}^{\pm}\dirac^{\pm}.
\end{eqnarray}

Combining Proposition~\ref{diracmh} with operator relation
(\ref{DV}), we have for $M_{k-s}^{\pm} \in {\cal M}_{k-s}^{\pm}$,
\begin{eqnarray}
\label{DDf}(\dirac^{\pm})^2
((mh)^2M_{k-s}^{\pm})&=&\dirac^{\pm}\left(-2V_{h,n/2}^{\pm}((mh)M_{k-s}^{\pm})+2(mh)V_{h,n/2}^{\pm}M_{k-s}^{\pm}\right)
\nonumber \\
&=& (-2)^2\left(V_{h,n/2+1}^{\pm}V_{h,n/2}^{\pm}M_{k-s}^{\pm}-V_{h,n/2}^{\pm}V_{h,n/2}^{\pm}M_{k-s}^{\pm}\right) \nonumber \\
&=& (-2)^2V_{h,n/2}^{\pm}M_{k-s}^{\pm} .
\end{eqnarray}
and
\begin{eqnarray}
\label{DDDf}(\dirac^{\pm})^3 ((mh)^3M_{k-s}^{\pm}) \nonumber
&=&(\dirac^{\pm})^2\left(-2V_{h,n/2}^{\pm}((mh)^2M_{k-s}^{\pm})\right. \nonumber \\
&& \left. +2(mh)V_{h,n/2}^{\pm}((mh)^2M_{k-s}^{\pm})  \right) \nonumber \\
&=&(-2)^3V_{h,n/2+2}^{\pm}V_{h,n/2}^{\pm}M_{k-s}^{\pm}\nonumber \\
&&+2(\dirac^{\pm})^2\left((mh)V_{h,n/2}^{\pm}((mh)^2M_{k-s}^{\pm}))\right)
\nonumber \\
&=&
(-2)^3\left(V_{h,n/2+2}^{\pm}V_{h,n/2}^{\pm}-V_{h,n/2+1}^{\pm}V_{h,n/2}^{\pm}\right.
\nonumber \\
&& \left. +V_{h,n/2}^{\pm}V_{h,n/2}^{\pm}\right)M_{k-s}^{\pm} \nonumber \\
&=& (-2)^3V_{h,n/2+1}^{\pm}V_{h,n/2}^{\pm}M_{k-s}^{\pm}
\end{eqnarray}

Continuing this procedure, we obtain by recursion
\begin{eqnarray} \label{Dsf}(\dirac^{\pm})^{s}
((mh)^sM_{k-s}^{\pm})&=& (-2)^{s}
V_{h,n/2+s-2}^{\pm}V_{h,n/2+s-3}^{\pm} \ldots V_{h,n/2}^{\pm}
M_{k-s}^{\pm}.
\end{eqnarray}

From this follows also $(mh)^sM_{k-s}^{\pm} \in \ker
(\dirac^{\pm})^{s+1}$. Formula~(\ref{Dsf}) gives us a motivation to
find explicit formulae for the polynomials $M_{k-s}^{\pm}$.  To this
end we need an explicit formula for the inverse of the iterated
composite operator $V_{h,n/2+s-2}^{\pm}V_{h,n/2+s-3}^{\pm} \ldots
V_{h,n/2}^{\pm}$. This means that we have to find an explicit
formula for the inverse of the operator $V_{h,r}^{\pm}$.
Unfortunately, we are only able to get an explicit formula for the
operator $R^{\pm}_{h,r}$ (c.f.~\cite{Ren}).

\begin{theorem}
For a lattice function $f_h:\BR^n_h \rightarrow \cl_{0,n}$ and for
$r>0$, the difference operator $J^{\pm}_{h,r}$ defined by $$
(J^{\pm}_{h,r}f_h)(mh) = \sum_{th \in [0,1]_h^\pm}
hd_h^{\pm}\left((th \mp h)^{(r-1)}_{\mp}f_h((th)(mh))\right)$$
satisfies
$$J^{\pm}_{h,r}R^{\pm}_{h,r}=I=R^{\pm}_{h,r}J^{\pm}_{h,r}.$$
Hereby we denote $[0,1]_h^+=[0,1)_h, [0,1]_h^-=(0,1]_h,$ and
$$(d_h^{\pm}g)(th):=\mp \frac{g(th)-g(th \pm h)}{h}.$$
\end{theorem}

\proof{\bf (c.f. \cite{Ren})} For $f_h:\BR^n_h \rightarrow
\cl_{0,n}$ and $r >0$,
$$ f_h(mh)=\sum_{th \in [0,1]_h}
hd_h^{\pm}\left((th \mp h)^{(r)}_{\mp}f_h((th)(mh))\right) $$

By a direct calculation,
\begin{eqnarray*}
&&d_h^{\pm}\left((th \mp h)^{(r)}_{\mp}f_h((th)(mh))\right)\\
&=&r(th \mp h)^{(r-1)}_{\mp}f_h((th)(mh)) +(th )^{(r)}_{\mp}(d_h^{\pm}f_h)((th)(mh))  \\
&=&(th \mp
h)^{(r-1)}_{\mp}\left(rf_h((th)(mh))+th(d_h^{\pm}f_h)((th)(mh))\right)
\end{eqnarray*}

On the other hand, applying the difference version of the chain
rule and the relation $\sum_{i=1}^n(m_ih) \dif^{\pm
i}=\euler-A_h^{\pm}$, we obtain

\begin{eqnarray*}
th(d_h^{\pm}f_h)((th)(mh))&=&\sum_{i=1}^n (th)(m_ih)
(\partial_{th}^{\pm
i}f_h)((th)(mh)) \\
&=& (E_{th}^{\pm}f_h)((th)(mh))-(A_{th}^{\pm}f_h)((th)(mh))
\end{eqnarray*}

Therefore,
\begin{eqnarray*}
f_h(mh) & = &
r(J^{\pm}_{h,r}f_h)(mh)+(E_{h}^{\pm}J^{\pm}_{h,r}f_h)(mh)-(A_{h}^{\pm}J^{\pm}_{h,r}f_h)(mh)\\
& = & (R^{\pm}_{h,r}J^{\pm}_{h,r}f_h)(mh).
\end{eqnarray*}

From the above two identities and the definitions of
$R^{\pm}_{h,r}$ and $J^{\pm}_{h,r}$, we get
$$J^{\pm}_{h,r}R^{\pm}_{h,r}=I=R^{\pm}_{h,r}J^{\pm}_{h,r}.$$
 \qed

Now, the construction of the inverse for $V^{\pm}_{h,r}$ seems to
be obvious. But, the obvious choice
$$(W^{\pm}_{h,r}f_h)(mh) = \sum_{th \in [0,1]_h^\pm} hd_h^{\pm}\left((th)^{(r-1)}_{\mp}f_h((th)(mh))\right)$$
is not an inverse of $V^{\pm}_{h,r}$, which can be easily checked
in the following way.

If we use the same technique as above, we obtain
$$f_h(mh) = \sum_{th \in [0,1]_h^\pm} hd_h^{\pm}\left((th)^{(r-1)}_{\mp}f_h((th)(mh))\right)$$
and by direct calculation
\begin{eqnarray*}
&&d_h^{\pm}\left((th)^{(r)}_{\mp}f_h((th)(mh))\right)\\
&=&r(th)^{(r-1)}_{\mp}f_h(thx) +(th \pm h)^{(r)}_{\mp}(d_h^{\pm}f_h)((th)(mh))   \\
&=&(th)^{(r-1)}_{\mp}\left(rf_h((th)(mh))+(th \pm h)
(d_h^{\pm}f_h)((th)(mh))\right)
\end{eqnarray*}

On the other hand, applying the difference version of the chain
rule and the relation $\sum_{i=1}^n (m_ih) \dif^{\pm
i}=\euler-A_h^{\pm}$, we obtain
\begin{eqnarray*}
(th \pm h)(d_h^{\pm}f_h)((th)(mh))&=&\sum_{i=1}^n ((th) (m_ih) \pm
hm_i) (\partial_{th}^{\pm
i}f_h)((th)(mh)) \\
&=& (E_{th}^{\pm}f_h)((th)(mh))-(A_{th}^{\pm}f_h)((th)(mh)) \\
&&\pm h \sum_{i=1}^n m_ih(\partial_{th}^{\pm i}f_h)((th)(mh)),
\end{eqnarray*}
but
\begin{eqnarray*}
\pm h \sum_{i=1}^n (m_ih)(\partial_{th}^{\pm i}f_h)((th)(mh))
&\neq & \pm h \sum_{i=1}^n (\partial_{th}^{\pm i}f_h)((th)(mh))\\
&& =(B^{\pm}_{th}f_h)((th)(mh)).
\end{eqnarray*}

\subsection{Difference operator calculus} \label{diffcalc}

Now we will establish some properties for our difference operators
introduced in Section~\ref{EulerGamma}.

Using the difference properties
$$(m_ih)\dif^{\pm i}(mh \mp h \e_i)_{\mp}^{(\alpha)}=\alpha_i(mh)_{\mp}^{(\alpha)}$$
and
$$(m_ih)\dif^{\pm i}\dif^{\mp i}(mh \mp h \e_i)_{\mp}^{(\alpha)}=(\alpha_i-1)(m_ih)\dif^{\pm i}(mh \mp h \e_i)_{\mp}^{(\alpha)}$$
we obtain by direct calculation the following formulae for
homogeneous polynomials of degree $k$, $P_k^{\pm} \in \Pi_k^{\pm}$,
\begin{eqnarray}
\label{BP}B_h^{\pm}P_k^{\pm} & = & \pm khP_k^{\pm}, \\
\label{AP}A_h^{\pm}P_k^{\pm} & = & \frac{kh^2}{1 \pm h}P_k^{\pm}, \\
\label{RP}R_{h,r}^{\pm}P_k^{\pm} & = & \left(r+k-\frac{kh^2}{1 \pm h} \right)P_k^{\pm}, \\
\label{VP}V_{h,r}^{\pm}P_k^{\pm} & = & \left(r+\left(1\pm
\frac{h}{2}\right)k-\frac{kh^2}{1 \pm h} \right)P_k^{\pm}.
\end{eqnarray}
Now, using the difference rules (\ref{difRule1}) and
(\ref{difRule2}), we get
\begin{eqnarray}
\label{Bmh}B_h^{\pm}((mh)f_h)&=&(mh)B_h^{\pm}f_h + h {\bf
1}^{\pm}f_h(mh)+h^2
\dirac^{\pm}f_h, \\ \nonumber \\
\label{Cmh}C_h^{\pm}((mh)f_h)&=&(mh)C_h^{\pm}f_h-\euler f_h, \\ \nonumber \\
\label{Amh}A_h^{\pm}((mh)f_h)&=&(mh)A_h^{\pm}f_h \mp h C_h^{\pm}f_h , \\ \nonumber \\
\label{Emh}\euler((mh)f_h)&=&(mh)\euler f_h+C_h^{\pm}f_h, \\ \nonumber \\
\label{Rmh}R_{h,r}^{\pm}((mh)f_h)&=&(mh)R_{h,r}^{\pm}f_h+(1 \pm
h)C_h^{\pm}f_h,
\\ \nonumber \\ \label{Vmh}V_{h,r}^{\pm}((mh)f_h)&=&(mh)V_{h,r}^{\pm}f_h+(1 \pm
h)C_h^{\pm}f_h + \nonumber \\
&+&\frac{1}{2}( h {\bf 1}^{\pm}f_h(mh)+h^2 \dirac^{\pm}f_h),
\end{eqnarray}
where ${\bf 1}^{\pm}:=\pm \sum_{i=1}^n \e_i.$ \\

As a direct consequence of formulae (\ref{VP}) and (\ref{Emh}), we
obtain for the discrete homogeneous monogenic polynomials of degree
$k$, $M_{k}^{\pm} \in {\cal M}_{k}^{\pm}$,
\begin{eqnarray}
\label{Gamh}\gam((mh)M_{k}^{\pm})&=&-\euler((mh)M_{k}^{\pm})-(mh)\dirac^{\pm}((mh)M_{k}^{\pm})
\nonumber \\
&=&\left(n+k-\frac{2kh^2}{1 \pm h}\pm
hk\right)(mh)M_k^{\pm}-C_h^{\pm}M_k^{\pm}
\end{eqnarray}
by applying Theorem~\ref{diracmh} and relation~(\ref{id2}).

Moreover, using identities (\ref{Dmh}),
(\ref{BP}),(\ref{RP}),(\ref{VP}),(\ref{id2}) and
Proposition~\ref{diracmh}, we obtain the relation
\begin{eqnarray} \label{Dhmh}\left(\dirac^{\pm}
(mh)\right)((mh)M_k^{\pm})&=&(-2R_{h,n/2}^{\pm}+\euler+\gam)
((mh)M_k^{\pm}) \nonumber \\
&=& \pm hk(mh)M_k^{\pm}-(2\pm 2h)C_h^{\pm}M_k^{\pm}.
\end{eqnarray}

Applying relations (\ref{Emh}) and (\ref{Cmh}) we have
\begin{equation}\label{eulermhmh} \euler((mh)^2M_k^{\pm})=k(mh)^2M_k^{\pm}+2(mh)C_h^{\pm}M_k^{\pm}-kM_k^{\pm}.  \end{equation}

From Theorem \ref{diracmh} and formulae
(\ref{Bmh}) and (\ref{Dhmh}), we get
\begin{equation}\label{Dirmhmh} \dirac^{\pm}((mh)^2M_k^{\pm})=-(2\pm 2h)C_h^{\pm}M_k^{\pm}+h {\bf
1}^{\pm}M_k^{\pm}.  \end{equation}

Using (\ref{id2}) and formulae (\ref{eulermhmh}) and
(\ref{Dirmhmh}), we obtain
\begin{equation}\label{gammhmh} \gam((mh)^2M_k^{\pm})=-k(mh)^2M_k^{\pm}\pm 2(mh)C_h^{\pm}M_k^{\pm}+kM_k^{\pm}-h {\bf
1}^{\pm}(mh)M_k^{\pm}.  \end{equation}

Applying recursively our formulae (\ref{Bmh})-(\ref{Vmh}), it is
possible to obtain explicit formulae for
$\dirac^{\pm}((mh)^sM_k^{\pm})$, $\euler((mh)^sM_k^{\pm})$ and
$\gam((mh)^sM_k^{\pm})$ by induction.

\subsection{Homogeneous powers} \label{Homopowers}

Contrary to the continuous case, the classical product between the
variable $mh$ and the homogeneous polynomial $P^{\pm}_k$ is not
homogeneous. However, applying formulae (\ref{Stirl1}) and
(\ref{Stirl2}) proved in Theorem~\ref{comb2}, we can say that the
product $(mh)P^{\pm}_k$ can be expressed as a linear combination of
homogeneous polynomials up to degree $k+1$. On the other hand, the
powers $x^s=(mh)^s$ are not homogeneous. For this purpose, we will
introduce the discrete analogues of $x^s$ in the following way:

Starting from the multi-index factorial powers of degree $|\alpha|$,
$(mh)^{(\alpha)}_{\mp}$, we introduce the polynomials $H^{\pm}_s, s
\in \BN$ by
\begin{equation}
\label{Homo}H^{\pm}_s(mh)=\left\{
\begin{array}{ccc}
\sum_{|\alpha|=k}\frac{(-1)^k k!}{\alpha!}(mh)^{(2 \alpha)}_{\mp}&
~if~s=2k, k \in \BN_0\\ \
\\\sum_{|\alpha|=k}\sum_{i=1}^n\frac{(-1)^k k!}{
\alpha!}(mh)^{(2\alpha + \e_i)}_{\mp} \e_i & ~if~s=2k+1, k \in
\BN_0.
\end{array}
\right.
\end{equation}

As a direct consequence of the identity
$$
x^s=\left\{
\begin{array}{ccc}
(-1)^k |x|^{2k}  & ~if~s=2k, k \in \BN_0\\
xx^{2k} &  ~if~s=2k+1, k \in \BN_0
\end{array}
\right.
$$
we can conclude by Lemma~\ref{prop3} and by the multinomial theorem,
that the polynomials $H_s^{\pm}$ give rise to homogeneous
polynomials of degree $s$ which approximate the powers $x^s=(mh)^s$
for small mesh width $h >0$.

As a direct consequence, the operator formulae proved in Subsection
\ref{diffcalc} are fulfilled for the powers $H^{\pm}_s$ and,
moreover, by direct calculation, we obtain the additional properties
$$ C^{\pm}_hH^{\pm}_s=H^{\pm}_{s+1},$$
and
$$\dirac^{\pm}H^{\pm}_s=-sH^{\pm}_{s-1}.$$
Let us remark that the term $C^{\pm}_hH^{\pm}_s$ is the discrete
version of the multiplication $xx^s$ in the continuous case.

\section{A discrete harmonic Fischer decomposition}

According to the classical theory of the finite differences, the
usual approximation of the Laplacian is given by
\begin{eqnarray}
\label{difLap} (\lapl u) (mh)&=&\sum_{i=1}^n
\frac{u(mh+h\e_i)+u(mh-h\e_i)-2u(mh)}{h^2} \nonumber \\
&=&\sum_{i=1}^n (\dif^{\mp i}\dif^{\pm i}u)(mh).
\end{eqnarray}

The first problem that arises now is that not all of our partial
difference operators do commute in the certain sense
(c.f.~\cite{GS,GS2}) and, moreover, we have no factorization of the
discrete Laplacian $\lapl$ by means of our difference Dirac operators
considered above, that is $\dirac^{\mp}\dirac^{\pm} \neq
-\e_0\lapl$.

Let us restrict ourselves in this section to the case of
quaternion-valued functions defined on lattices in $\BR^3$. \\
Let us remark that the quaternionic variable $mh$ is identified with
the $4 \times 4$ matrix
$$
mh=\left(
\begin{array}{cccc}
0 & -m_1h& -m_2h & -m_3h \\
m_1h     & 0    &   -m_3h & m_2h      \\
m_2h     & m_3h & 0      &  m_1h    \\
m_3h     & -m_2h& m_1h      & 0
\end{array}
\right).
$$

In~\cite{GH} for a lattice function $f_h: \BR^3_h \rightarrow
\mathbb{H}$ given by $$f_h=\sum_{i=0}^3 f_h^i \e_i=f_h^0\e_0+\Vec
f_h$$ a finite difference approximation of our Dirac operator was
defined in the form
\begin{eqnarray}
\label{dimp}\dirac^{-+}f_h&=& \left(
\begin{array}{cccc}
0 & -\dif^{-1}& -\dif^{-2} & -\dif^{-3} \\
\dif^{-1}     & 0    &   -\dif^{3} & \dif^{2}      \\
\dif^{-2}    & \dif^{3} & 0      &  -\dif^{1}    \\
\dif^{-3}     & -\dif^{2}& \dif^{1}      & 0
\end{array}
\right) \left(\begin{array}{c}
f_h^0  \\
f_h^1        \\
f_h^2       \\
f_h^3
\end{array}
\right) \nonumber \\&=& \left(
\begin{array}{c}
-\Div_h^- \Vec f_h\\
\Grad_h^- f_h^0+\Curl_h^+\Vec f_h
\end{array}
\right) \\ \nonumber \\ \nonumber \\ \label{dipm}\dirac^{+-}f_h
&=&\left(
\begin{array}{cccc}
0 & -\dif^{1}& -\dif^{2} & -\dif^{3} \\
\dif^{1}     & 0    &   -\dif^{-3} & \dif^{-2}      \\
\dif^{2}    & \dif^{-3} & 0      &  -\dif^{-1}    \\
\dif^{3}     & -\dif^{-2}& \dif^{-1}      & 0
\end{array}
\right)  \left(\begin{array}{c}
f_h^0  \\
f_h^1        \\
f_h^2       \\
f_h^3
\end{array}
\right)\nonumber \\ &=&\left(
\begin{array}{c}
-\Div_h^+ \Vec f_h\\
\Grad_h^+ f_h^0+\Curl_h^-\Vec f_h
\end{array}
\right)
\end{eqnarray}
with $\Div_h^{\pm}\Vec f_h=\sum_{i=1}^3\dif^{\pm i}f_h^i$,
$\Grad_h^{\pm} f_h^0=\sum_{i=1}^3(\dif^{\pm i}f_h^0)\e_i$ and
$$
\Curl_h^{\pm}\Vec f_h=\left|
\begin{array}{ccc}
\e_1 & \e_2 & \e_3\\
\dif^{\pm 1} & \dif^{\pm 2}& \dif^{\pm 3}\\
f_h^1 & f_h^2 & f_h^3
\end{array}
\right|.
$$

In the latter form one can easily see the similarity with the usual
Dirac operator
$$
D f=\left(
\begin{array}{c}
-\Div \Vec f\\
\Grad \Sc f+\Curl \Vec f
\end{array}
\right).
$$

Using the discrete identities
\begin{eqnarray*}
\Div_h^{\pm}\Curl_h^{\pm} \Vec f_h=0 \\
\Curl_h^{\pm}\Grad_h^{\pm} f_h^0={\bf 0} \\
\Curl_h^{\pm}\Curl_h^{\mp} \Vec f_h=-\lapl \Vec f_h+\Grad_h^{\mp}
\Div_h^{\pm} \Vec f_h
\end{eqnarray*}
we obtain the following factorization of the discrete Laplacian
\begin{eqnarray}
\label{factlapl}\dirac^{+-}\dirac^{-+}f_h=\left(
\begin{array}{c}
-\lapl f_h^0\\
-\lapl \Vec f_h
\end{array}
\right)=\dirac^{-+}\dirac^{+-}f_h.
\end{eqnarray}

Now, we are able to obtain a Fischer decomposition for the discrete
Dirac operators $\dirac^{-+}$ and $\dirac^{+-}$.

Using the fact that
\begin{eqnarray*}\dirac^{-+}(mh)_{+}^{(\alpha)}&=&\left(
\begin{array}{c}
 0\\
\Grad_h^- (mh)_{+}^{(\alpha)}
\end{array}
\right)\\&=&\left(
\begin{array}{c}
 0\\
\alpha_1(mh)_{+}^{(\alpha-\e_1)} \\
\alpha_2(mh)_{+}^{(\alpha-\e_2)} \\
\alpha_3(mh)_{+}^{(\alpha-\e_3)}
\end{array}
\right)
\end{eqnarray*}
as well as
\begin{eqnarray*}\dirac^{+-}(mh)_{-}^{(\alpha)}&=&\left(
\begin{array}{c}
 0\\
\Grad_h^+ (mh)_{-}^{(\alpha)}
\end{array}
\right)\\&=&\left(
\begin{array}{c}
 0\\
\alpha_1(mh)_{-}^{(\alpha-\e_1)} \\
\alpha_2(mh)_{-}^{(\alpha-\e_2)} \\
\alpha_3(mh)_{-}^{(\alpha-\e_3)}
\end{array}
\right)
\end{eqnarray*}
we can prove as in Lemma ~\ref{prop5}, the inclusion properties
$\dirac^{+-} \Pi_k^{+} \subset \Pi_{k-1}^{+}$, $ \dirac^{-+}
\Pi_k^{-} \subset \Pi_{k-1}^{-}$ and, moreover, replacing
$\dirac^{+-}$ by $\dirac^{+}$ and $\dirac^{-+}$ by $\dirac^{-}$ in
the inner product~(\ref{innerfischer}), we obtain the Fischer
decompositions:

\begin{theorem}{\bf Fischer decomposition for $\dirac^{-+}$ and $\dirac^{+-}$}\label{fisherD}\\ \ \\ Let $P^{-}_k \in \Pi_{k}^{-}$ (respectively, $P^{+}_k \in \Pi_{k}^{+}$)  then
\begin{eqnarray}
P^{-}_k&=&\sum_{s=0}^{k-1}(mh)^sM_{k-s}^{-+}, \\
P^{+}_k&=&\sum_{s=0}^{k-1}(mh)^sM_{k-s}^{+-}.
\end{eqnarray}
where each $M^{-+}_j$( respectively, $M^{+-}_j$) denotes a
homogeneous discrete monogenic polynomial of degree $j$, that is,
$M^{-+}_j \in \Pi_{j}^{-} \cap \ker \dirac^{-+}$(respectively,
$M^{+-}_j \in \Pi_{j}^{+} \cap \ker \dirac^{+-}$).
\end{theorem}

From the factorization property (\ref{factlapl}), we have
$$[(mh)^2 P_{k}^{\pm},Q_{k}^{\pm}]_h=-[ P_{k}^{\pm},\lapl
Q_k^{\pm}]_h,$$ which allows us to obtain the Fischer decomposition
for the discrete Laplacian:

\begin{theorem}{\bf Fischer decomposition for $\lapl$}\label{fischerlapl}\\ \ \\
Let $P^{\pm}_k \in \Pi_{k}^{\pm}$   then
$$P^{\pm}_k=\sum_{2 s \leq k}|mh|^{2s}\mathcal{H}_{k-2s}^{\pm},$$
where each $\mathcal{H}_{j}^{\pm}$ denotes a homogeneous discrete
harmonic polynomial of degree $j$, that is, $\mathcal{H}_{j}^{\pm}
\in \Pi_{j}^{\pm} \cap \ker \lapl$.
\end{theorem}

As a consequence of Theorem~\ref{fisherD}, we obtain Fischer
decompositions which relate the discrete harmonic and the discrete
monogenic polynomials.

\begin{corollary}{\bf Fischer decomposition }
Let $\mathcal{H}_{k}^{\pm} \in \Pi_{k}^{\pm} \cap \ker \lapl$ then
\begin{eqnarray}
\mathcal{H}_{k}^{-}&=&M_{k}^{-+}+(mh)M_{k-1}^{-+}, \\
\mathcal{H}_{k}^{+}&=&M_{k}^{+-}+(mh)M_{k-1}^{+-}.
\end{eqnarray}
where each $M^{-+}_j$( respectively, $M^{+-}_j$) denotes a
homogeneous discrete monogenic polynomial of degree $j$, that is,
$M^{-+}_j \in \Pi_{j}^{-} \cap \ker \dirac^{-+}$(respectively,
$M^{+-}_j \in \Pi_{j}^{+} \cap \ker \dirac^{+-}$).
\end{corollary}

To define the Euler and Gamma operator $E_{h}^{+-},\Gamma_h^{+-}$
and $E_{h}^{-+},\Gamma_h^{-+}$ for the modified Dirac operators
$\dirac^{+-}$ and $\dirac^{-+}$, respectively, we start to
calculate the products $(mh)\dirac^{+-}f_h$ and
$(mh)\dirac^{+-}f_h$. By straightforward calculations we obtain
\begin{eqnarray} \label{idD-+} (mh)\dirac^{-+}f_h& = & -\left(
\begin{array}{cccc}
\dif^{-1}f_h^0& \dif^{-2}f_h^0 & \dif^{-3}f_h^0 \\
\dif^{-1}f_h^1&\dif^{2}f_h^1&\dif^{3}f_h^1      \\
\dif^{1}f_h^2&\dif^{-2}f_h^2&\dif^{3}f_h^2 \\
\dif^{1}f_h^3&\dif^{2}f_h^3&\dif^{-3}f_h^3
\end{array}
 \right)
 \left(
\begin{array}{c}
m_1h\\
m_2h\\
m_3h
\end{array}
 \right) \nonumber \\
 & + &\left|
\begin{array}{ccc}
\e_1 & \e_2 & \e_3\\
m_1h & m_2h& m_3h\\
\dif^{-1} & \dif^{-2} & \dif^{-3}
\end{array}
\right|f_h^0+\left|
\begin{array}{ccc}
\e_0 & \e_2 & \e_3\\
m_1h & m_3h& m_2h\\
\dif^{-1} & \dif^{3} & \dif^{2}
\end{array}
\right|f_h^1 \nonumber \\ & + & \left|
\begin{array}{ccc}
\e_0 & \e_1 & \e_3\\
-m_2h & -m_3h& m_1h\\
-\dif^{-2} & -\dif^{3} & \dif^{1}
\end{array}
\right|f_h^2 \nonumber \\ & + & \left|
\begin{array}{ccc}
\e_0 & \e_1 & \e_3\\
m_3h & m_2h& m_1h\\
\dif^{-3} & \dif^{2} & \dif^{1}
\end{array}
\right|f_h^3.
\end{eqnarray}
and
\begin{eqnarray}
\label{idD+-} (mh)\dirac^{+-}f_h & = & -\left(
\begin{array}{cccc}
\dif^{1}f_h^0&\dif^{2}f_h^0&\dif^{3}f_h^0 \\
\dif^{1}f_h^1&\dif^{-2}f_h^1&\dif^{-3}f_h^1      \\
\dif^{-1}f_h^2&\dif^{2}f_h^2&\dif^{-3}f_h^2 \\
\dif^{-1}f_h^3&\dif^{-2}f_h^3&\dif^{3} f_h^3
\end{array}
 \right)\left(
\begin{array}{c}
m_1h\\
m_2h\\
m_3h
\end{array}
 \right) \nonumber \\ & + &
 \left|
\begin{array}{ccc}
\e_1 & \e_2 & \e_3\\
m_1h & m_2h& m_3h\\
\dif^{1} & \dif^{2} & \dif^{3}
\end{array}
\right|f_h^0+\left|
\begin{array}{ccc}
\e_0 & \e_2 & \e_3\\
m_1h & m_3h& m_2h\\
\dif^{1} & \dif^{-3} & \dif^{-2}
\end{array}
\right|f_h^1 \nonumber \\ &+& \left|
\begin{array}{ccc}
\e_0 & \e_1 & \e_3\\
-m_2h & -m_3h& m_1h\\
-\dif^{2} & -\dif^{-3} & \dif^{-1}
\end{array}
\right|f_h^2 \nonumber \\ & + & \left|
\begin{array}{ccc}
\e_0 & \e_1 & \e_2\\
m_3h & m_2h& m_1h\\
\dif^{3} & \dif^{-2} & \dif^{-1}
\end{array}
\right|f_h^3.
\end{eqnarray}

Hence, we can define the difference Euler operators $E_h^{-+}$ and
$E_h^{+-}$ as \begin{eqnarray*} && (E_h^{-+}f_h)(mh)=\\
&&\left(
\begin{array}{cccc}
(\dif^{-1}f_h^0)(mh+h\e_1)&(\dif^{-2}f_h^0)(mh+h\e_2)&(\dif^{-3}f_h^0)(mh+h\e_3) \\
(\dif^{-1}f_h^1)(mh+h\e_1)&(\dif^{2}f_h^1)(mh-h\e_2)&(\dif^{3}f_h^1)(mh-h\e_3)      \\
(\dif^{1}f_h^2)(mh-h\e_1)&(\dif^{-2}f_h^2)(mh+h\e_2)&(\dif^{3}f_h^2)(mh-h\e_3) \\
(\dif^{1}f_h^3)(mh-h\e_1)&(\dif^{2}f_h^3)(mh-h\e_2)&(\dif^{-3}f_h^3)(mh+h\e_3)
\end{array}
 \right)\left(
\begin{array}{c}
m_1h\\
m_2h\\
m_3h
\end{array}
 \right)
 \end{eqnarray*}
and
\begin{eqnarray*}
&&(E_h^{+-}f_h)(mh)=\\
&&\left(
\begin{array}{cccc}
(\dif^{1}f_h^0)(mh-h\e_1)&(\dif^{2}f_h^0)(mh-h\e_2)&(\dif^{3}f_h^0)(mh-h\e_3) \\
(\dif^{1}f_h^1)(mh-h\e_1)&(\dif^{-2}f_h^1)(mh+h\e_2)&(\dif^{-3}f_h^1)(mh+h\e_3)      \\
(\dif^{-1}f_h^2)(mh+h\e_1)&(\dif^{2}f_h^2)(mh-h\e_2)&(\dif^{-3}f_h^2)(mh+h\e_3) \\
(\dif^{-1}f_h^3)(mh+h\e_1)&(\dif^{-2}f_h^3)(mh+h\e_2)&(\dif^{3}f_h^3)(mh-h\e_3)
\end{array}
 \right)\left(
\begin{array}{c}
m_1h\\
m_2h\\
m_3h
\end{array}
 \right)
 \end{eqnarray*}

The difference Gamma operators $\Gamma_h^{-+}$ and $\Gamma_h^{+-}$ are
defined by
\begin{eqnarray*}
&&(\Gamma_h^{-+}f_h)(mh)= \\
&&=h\left(
\begin{array}{cccc}
(\dif^{-1}\dif^{1}f_h^0)(mh) &(\dif^{-2}\dif^{2}f_h^0)(mh)&(\dif^{-3}\dif^{3}f_h^0)(mh)\\
(\dif^{-1}\dif^{1}f_h^1)(mh) &-(\dif^{-2}\dif^{2}f_h^1)(mh)&-(\dif^{-3}\dif^{3}f_h^1)(mh)    \\
-(\dif^{-1}\dif^{1}f_h^2)(mh) &(\dif^{-2}\dif^{2}f_h^2)(mh)&-(\dif^{-3}\dif^{3}f_h^2)(mh) \\
-(\dif^{-1}\dif^{1}f_h^3)(mh)
&-(\dif^{-2}\dif^{2}f_h^3)(mh)&(\dif^{-3}\dif^{3}f_h^3)(mh)
\end{array}
 \right)\left(
\begin{array}{c}
m_1h\\
m_2h\\
m_3h
\end{array}
 \right) \\ && - \left|
\begin{array}{ccc}
\e_1 & \e_2 & \e_3\\
m_1h & m_2h& m_3h\\
\dif^{-1} & \dif^{-2} & \dif^{-3}
\end{array}
\right|f_h^0(mh)-\left|
\begin{array}{ccc}
\e_0 & \e_2 & \e_3\\
m_1h & m_3h& m_2h\\
\dif^{-1} & \dif^{3} & \dif^{2}
\end{array}
\right|f_h^1(mh) \\ &&- \left|
\begin{array}{ccc}
\e_0 & \e_1 & \e_3\\
-m_2h & -m_3h& m_1h\\
-\dif^{-2} & -\dif^{3} & \dif^{1}
\end{array}
\right|f_h^2(mh) \\  && -  \left|
\begin{array}{ccc}
\e_0 & \e_1 & \e_2\\
m_3h & m_2h& m_1h\\
\dif^{-3} & \dif^{2} & \dif^{1}
\end{array}
\right|f_h^3(mh)
\end{eqnarray*}
and
 \begin{eqnarray*}&&(\Gamma_h^{+-}f_h)(mh)= \\ && = -h\left(
\begin{array}{cccc}
(\dif^{-1}\dif^{1}f_h^0)(mh) &(\dif^{-2}\dif^{2}f_h^0)(mh)&(\dif^{-3}\dif^{3}f_h^0)(mh)\\
(\dif^{-1}\dif^{1}f_h^1)(mh) &-(\dif^{-2}\dif^{2}f_h^1)(mh)&-(\dif^{-3}\dif^{3}f_h^1)(mh)    \\
-(\dif^{-1}\dif^{1}f_h^2)(mh) &(\dif^{-2}\dif^{2}f_h^2)(mh)&-(\dif^{-3}\dif^{3}f_h^2)(mh) \\
-(\dif^{-1}\dif^{1}f_h^3)(mh)
&-(\dif^{-2}\dif^{2}f_h^3)(mh)&(\dif^{-3}\dif^{3}f_h^3)(mh)
\end{array}
 \right)\left(
\begin{array}{c}
m_1h\\
m_2h\\
m_3h
\end{array}
 \right) \\ &&  \\ && -
 \left|
\begin{array}{ccc}
\e_1 & \e_2 & \e_3\\
m_1h & m_2h& m_3h\\
\dif^{1} & \dif^{2} & \dif^{3}
\end{array}
\right|f_h^0(mh)-\left|
\begin{array}{ccc}
\e_0 & \e_2 & \e_3\\
m_1h & m_3h& m_2h\\
\dif^{1} & \dif^{-3} & \dif^{-2}
\end{array}
\right|f_h^1(mh)\\ && - \left|
\begin{array}{ccc}
\e_0 & \e_1 & \e_3\\
-m_2h & -m_3h& m_1h\\
-\dif^{2} & -\dif^{-3} & \dif^{-1}
\end{array}
\right|f_h^2(mh) \\ &&  -  \left|
\begin{array}{ccc}
\e_0 & \e_1 & \e_2\\
m_3h & m_2h& m_1h\\
\dif^{3} & \dif^{-2} & \dif^{-1}
\end{array}
\right|f_h^3(mh).
\end{eqnarray*}

As in Subsection~\ref{EulerGamma}, we have
$(mh)\dirac^{-+}=-E_h^{-+}-\Gamma_h^{-+}$ (respectively,
$(mh)\dirac^{+-}=-E_h^{+-}-\Gamma_h^{+-}$) and the polynomials
$P_k^{\pm} \in \Pi_k^{\pm}$ satisfy $E_h^{-+}P_k^{-}=kP_k^{-}$,
(respectively, $E_h^{+-}P_k^{+}=kP_k^{+}$). Moreover, if $P_k^{-}
\in \ker \dirac^{-+}$ (respectively, $P_k^{+} \in \ker \dirac^{+-}$)
then we have $\Gamma_h^{-+}P_k^{-}=-kP_k^{-}$, (respectively,
$\Gamma_h^{+-}P_k^{+}=-kP_k^{+}$).

Like in Theorem~\ref{diracmh} we can prove the operator property
$\dirac^{-+}E^{-+}_h=I+E^{-+}_h\dirac^{-+}$ (respectively,
$\dirac^{+-}E^{+-}_h=I+E^{+-}_h\dirac^{+-}$). In the same way we get
analogous relations to the ones presented in
Subsection~\ref{EulerGamma} and in Subsection~\ref{diffcalc}. In
addition it is also possible define the discrete versions of the
quaternionic powers $(mh)^s$ with respect to our difference Dirac
operators $\dirac^{-+}$ and $\dirac^{+-}$, using a similar
construction as in Subsection \ref{Homopowers}.

At this point it would be interesting to know if the operator
setting we discussed here for the quaternionic case has an
equivalent operator setting in the general case of Clifford
algebras. Up to know it is not known, but we will discuss it in a
forthcoming paper~\cite{FKS}.

\end{document}